\documentclass[12pt,reqno]{amsart}

\usepackage{amssymb}
\usepackage[T1]{fontenc}
\usepackage[latin1]{inputenc} 
\usepackage{dsfont, mathrsfs}
\usepackage{stmaryrd}
\usepackage{amsthm}
\usepackage{graphicx}
\usepackage{color, xcolor, soul}
\usepackage{xspace}
\usepackage{yhmath, epigraph}
\usepackage{caption, subcaption}
\usepackage{relsize}
\usepackage{shuffle}
\usepackage{float, adjustbox} 
\usepackage{mleftright}
\usepackage{pdfpages}
\usepackage{pax} 
\usepackage[
        colorlinks,
        hyperfootnotes=true,
        pagebackref, 
        hyperindex,
        linkcolor=blue]{hyperref}
\usepackage{geometry}
\usepackage{framed} 
\colorlet{shadecolor}{gray!40}
\usepackage{bbm}
\usepackage[bbgreekl]{mathbbol}
\DeclareSymbolFontAlphabet{\mathbb}{AMSb}
\DeclareSymbolFontAlphabet{\mathbbl}{bbold}
\usepackage{tikz}
\usetikzlibrary{positioning}
\usetikzlibrary{calc}
\usetikzlibrary{arrows.meta}
\usetikzlibrary{arrows}
\usepackage{pgfplots}
\pgfplotsset{compat=1.18}
\makeatletter 
\newcommand\HREF[2]{\hyper@linkurl{#2}{#1}}
\makeatother  
\usepackage{verbatim}
\DeclareFontFamily{U}{mathx}{\hyphenchar\font45}
\DeclareFontShape{U}{mathx}{m}{n}{
      <5> <6> <7> <8> <9> <10>
      <10.95> <12> <14.4> <17.28> <20.74> <24.88>
      mathx10
      }{}
\DeclareSymbolFont{mathx}{U}{mathx}{m}{n}
\DeclareFontSubstitution{U}{mathx}{m}{n}
\DeclareMathAccent{\widecheck}{0}{mathx}{"71}
\DeclareMathAccent{\wideparen}{0}{mathx}{"75}

\newcommand{\Zindex}[1]{}  
\usepackage{listings} 
\usepackage{marginnote}
\allowdisplaybreaks

\theoremstyle{plain}
\newtheorem{theorem}{Theorem}[section]
\newtheorem{lemma}[theorem]{Lemma}

\newtheorem{corollary}[theorem]{Corollary}

\theoremstyle{definition}

\theoremstyle{remark}
\newtheorem{remark}[theorem]{Remark}

\makeatletter
\newcommand*{\trimabs}[1]{%
  \mathpalette{\@trimleftright{|}{|}}{#1}%
}
\newcommand*{\@trimleftright}[4]{%
  \sbox0{$#3#4\m@th$}%
  \sbox2{$#3\vcenter{}$}
  \dimen0=\ht0 %
  \ifdim\dimen0<\ht2 %
    \dimen0=\ht2 %
  \fi
  \dimen2=\dp0 %
  \ifdim\dimen2<\z@
    \dimen2=\z@
  \fi
  \dimen0=\dimexpr(\dimen0-\dimen2)/2 -\ht2\relax
  \raisebox{\dimen0}{%
    $#3\mleft#1\raisebox{-\dimen0}{\box0}\mright#2\m@th$%
  }%
}
\makeatother

\newcommand{\Bb}{\mathbb{B}}
\newcommand{\Cc}{\mathbb{C}}

\newcommand{\Ee}{\mathbb{E}}

\newcommand{\Pp}{\mathbb{P}}

\newcommand{\Rr}{\mathbb{R}}

\newcommand{\Un}{\mathds{1}}

\newcommand{\Ce}{\mathcal{C}}
\newcommand{\De}{\mathcal{D}}
\newcommand{\Eee}{\mathcal{E}}
\newcommand{\Fe}{\mathcal{F}}
\newcommand{\Ge}{\mathcal{G}}

\newcommand{\Le}{\mathcal{L}}

\newcommand{\Ze}{\mathcal{Z}}

\newcommand{\Feb}{{{\boldsymbol{\mathcal{F}}}}}

\newcommand{\Zeb}{{{\boldsymbol{\mathcal{Z}}}}}

\newcommand{\Ns}{{\mathscr{N}}}

\newcommand{\Us}{{\mathscr{U}}}

\newcommand{\Zg}{{\mathfrak{Z}}}

\newcommand{\gammab}{{\boldsymbol{\gamma}}}

\newcommand{\mub}{{\boldsymbol{\mu}}}

\newcommand{\Gammab}{{\boldsymbol{\Gamma}}}

\newcommand{\Bd}{{\boldsymbol{B}}}

\newcommand{\Fb}{{\boldsymbol{F}}}

\newcommand{\Tb}{{\boldsymbol{T}}}

\newcommand{\Vb}{{\boldsymbol{V}}}
\newcommand{\Wb}{{\boldsymbol{W}}}

\newcommand{\Yb}{{\boldsymbol{Y}}}
\newcommand{\Zb}{{\boldsymbol{Z}}}

\newcommand{\vb}{{\boldsymbol{v}}}

\newcommand{\xb}{{\boldsymbol{x}}}

\newcommand{\zb}{{\boldsymbol{z}}}

\def\ee{ { \mathbbm{e} } }

\newcommand{\ensemble}[1]{ {\left\lbrace #1 \right\rbrace } } 
\newcommand{\prth}[1]{\!\left( #1 \right) }
\newcommand{\crochet}[1]{\!\left[ #1 \right] }  
\newcommand{\intcrochet}[1]{\llbracket #1 \rrbracket} 
\newcommand{\abs}[1]{{\left| #1 \right|}}

\newcommand{\Esp}[1]{ \Ee  \prth{ #1 } }  
\newcommand{\Prob}[1]{ \Pp \prth{ #1 } } 
 
\def\inv{^{-1}}

\newcommand{\cvlaw}[2]{\xrightarrow[#1 \, \rightarrow \, #2]{\Le}}

\newcommand{\equivalent}[1]{ {\underset{#1 }{\sim} } }

\newcommand{\bracket}[1]{\left\langle #1 \right\rangle}
\newcommand{\Unens}[1]{ \Un_{ \ensemble{#1} } }
\newcommand{\pe}[1]{{\left[ #1 \right]}}

\newcommand{\transp}[1]{\vphantom{#1^t}^t \!\;\! #1 } 

\def\eqlaw{\stackrel{\Le}{=}}

\newcommand{\centered}[1]{{\overset{\raisebox{-0.1cm}{\tiny$\circ$}}{#1}}}

\def\Arg{{ \operatorname{Arg} }}

\def\IID{{ \operatorname{i.i.d.} }}

\def\Exp{ \operatorname{Exp} }
\def\Ber{ \operatorname{Ber} }

\def\Gumbel{ \operatorname{Gb} }

\def\Argtanh{{\operatorname{Argtanh}}}

\def\geq{\geqslant}
\def\leq{\leqslant}

\let\oldforall\forall
\def\forall{\oldforall\,} 

\let\oldexists\exists
\def\exists{\oldexists\,}

\newcommand{\emailhref}[1]{ \email{\href{mailto:#1}{#1}} }

\title[The decoupling field approach to the Curie-Weiss model]{Decoupling the i.i.d. field and the randomisation \\ field in the Curie-Weiss model}
\author[Y. Barhoumi-Andr\'eani]{Yacine Barhoumi-Andr\'eani}
\emailhref{yacine.barhoumi@math.bas.bg}
\address{Bulgarian Academy of Sciences,       
Institute of Mathematics and Informatics, 
Department of Algebra and Logic, 
Acad. Georgi Bonchev Str., Block 8, 1113 Sofia, Bulgaria.}

\author[P. Eichelsbacher]{Peter Eichelsbacher}
\emailhref{peter.eichelsbacher@rub.de}
\address{Ruhr-Universit\"at Bochum, 
Fakult\"at f\"ur Mathematik, 
Universit\"atsstrasse 150, 
44780 Bochum, Germany.}
\date{\today}


\begin{document}
\begin{abstract}
Using the De Finetti representation of the Curie-Weiss model, the uniform coupling of Bernoulli random variables and the Laplace inversion formula (almost surely), we show that the full phase diagram of the Curie-Weiss model can be explained by a competition between the De Finetti randomisation and an approximate Gaussian process indexed by a complex variable that is equal to the inverse Laplace transform on a complex line of a Brownian Bridge. A more refined process type of rescaling shows that this is a modification of the Brownian Sheet that is at the core of all Gaussian random variables in the limits obtained in the model. This almost sure Laplace inversion approach allows moreover to treat all types of spin laws in the same vein as the Curie-Weiss Bernoulli spins.

This gives a natural explanation of several results that already appeared in the literature in the subcritical and critical case in addition to produce new analogous results in the super-critical case.

The functional approach here defined can moreover be extended to a wide class of statistical mechanical models that includes the Ising model in any dimension. 
\end{abstract}
\maketitle

\setcounter{tocdepth}{1}
\tableofcontents 

\section{Introduction}\label{Sec:Intro}

\subsection{Motivations and history}\label{SubSec:Intro:Motivations}

The Curie-Weiss model of $n$ spins at inverse temperature $ \beta \geq 0 $ (and external magnetic field $ \mu = 0 $) is the law of the random variables $ (X_k^{(\beta)})_{1 \leq k \leq n} $ defined by 
\begin{align}\label{Def:CurieWeiss}
\Pp_n^{(\beta)} \equiv \Pp_{ \prth{ X_1^{(\beta)}\!, \dots,\, X_n^{(\beta)} } } := \frac{e^{ \frac{\beta}{2n} S_n^2 } }{ \Esp{ e^{\frac{\beta}{2n} S_n^2 } } } \bullet \Pp_{ (X_1, \dots, X_n) }
\end{align}
where $ f \bullet \Pp $ denotes the bias/penalisation/tilting of $ \Pp $ by $f$ and, in the classical version of the model,
\begin{align*}
(X_k)_{1 \leq k \leq n}\sim \textrm{i.i.d. Ber}_{\ensemble{\pm 1}}(1/2), \qquad\qquad S_n := \sum_{k = 1}^n X_k
\end{align*}

Its \textit{unnormalised} magnetisation is classically given by
\begin{align}\label{Def:UnnormalisedMagnetisation}
M_n^{(\beta)} := \sum_{k = 1}^n X_k^{(\beta)}
\end{align}


Pierre Curie first proposed this simple statistical mechanical model in 1895 \cite{CurieThese} and Pierre-Ernest Weiss improved it in 1907 \cite{WeissOnCurie}. It is an exactly solvable model of ferromagnetism (ferromagnetic alloys can change their magnetic behaviour on their own when heated passed a critical temperature). 

A modern understanding of this model sees it as a \textit{mean-field} approximation of the more sophisticated Ising model (without external magnetic field) that uses an interaction between nearest neighbors $ \sum_{i \sim j} X_i X_j $, i.e.
\begin{align}\label{Def:Ising}
\widetilde{\Pp}_n^{(\beta)} \equiv \Pp_{ \prth{ \widetilde{X}_1^{(\beta)}\!, \dots,\, \widetilde{X}_n^{(\beta)} } } := \frac{e^{  \beta  \sum_{i \sim j} X_i X_j } }{ \Esp{ e^{ \beta  \sum_{i \sim j} X_i X_j } } } \bullet \Pp_{ (X_1, \dots, X_n) }
\end{align}

It can also be seen as an Ising model on the complete graph with $n$ vertices with a particular graph Laplacian. See for example Friedli and Velenik \cite[ch. 2]{FriedliVelenik} for a recent exposition or the more classical references \cite{Brout, KacMecaStat, StanleyMecaStat, Thompson}.


The links between the Curie-Weiss model \eqref{Def:CurieWeiss} and the Ising model \eqref{Def:Ising} are commonly presented in light of this mean field approximation paradigm; one more similarity will be presented in this paper: the existence of an \textit{i.i.d. field} and a \textit{randomisation field} (see Annex~\ref{Sec:Ising}).

\medskip

In view of \eqref{Def:CurieWeiss}, the Curie-Weiss spins are exchangeable. There is a well known De Finetti measure $ \widetilde{\nu}_{n, \beta} : \crochet{0, 1} \to \crochet{0, 1} $ that allows to write the spins as a mixture of $ \Ber_\ensemble{\pm 1}(p) $ random variables~:
\begin{align}\label{Eq:DeFinettiCW:Meas:P}
\Pp_n^{(\beta)} \equiv \Pp_{(X_1^{(\beta)}, \dots,\, X_n^{(\beta)} ) } = \int_{\crochet{0, 1}} \Pp_{ ( X_1(p), \dots,\, X_n(p) ) } \, \widetilde{\nu}_{n, \beta}(dp), \quad (X_k(p))_{1 \leq k \leq n} \sim \mathrm{i.i.d.} \Ber_{\ensemble{\pm 1}}(p)
\end{align}

The description of $ \widetilde{\nu}_{n, \beta} $ is explicitly given by (see the companion paper \cite[Lem.~B.1]{BarhoumiButzekEichelsbacher} for a proof or references cited)
\begin{align}\label{Eq:DeFinettiMeasureP}
\widetilde{\nu}_{n, \beta}(dp) = \widetilde{f}_{n, \beta}(p) dp, \qquad \widetilde{f}_{n, \beta}(p) := \frac{1}{\Ze_{n, \beta} }  e^{ -\frac{n}{2\beta} \Argtanh(2p - 1)^2 - (n/2 + 1) \ln(1 - (2p - 1)^2) }
\end{align}

The renormalisation constant $ \Ze_{n, \beta} $ is defined by the equality $ \int_0^1 \widetilde{f}_{n, \beta}(p) dp = 1 $ and we recall that $ \Argtanh(x) = \frac{1}{2} \log\abs{\frac{1 + x}{1 - x}} $ for $ \abs{x} < 1 $.

\medskip
\subsection{Main results}\label{SubSec:Intro:Result}

One aim of this paper is to introduce a new paradigm of computation to extract Gaussian and non-Gaussian information from the Curie-Weiss model. To achieve this goal, one first step consists in reproving the following theorem that is classical and that can be found e.g. in the books
\cite{Brout, FriedliVelenik, KacMecaStat, StanleyMecaStat, Thompson} or in the papers \cite{EllisNewman, EllisNewmanRosen} (the Bernoulli case is first due to Simon and Griffiths \cite[thm. 1]{SimonGriffiths}).

\begin{shaded}
\begin{theorem}[Complete phase transition of the Curie-Weiss model]\label{Theorem:MainIntro}
We have :
\begin{enumerate}

\medskip
\item If $ \beta < 1 $, 
\begin{align*}
\frac{1}{\sqrt{n}} \, M_n^{(\beta )} \cvlaw{n}{+\infty} \Ns\prth{ 0, \frac{1}{1 - \beta } }. 
\end{align*}

\medskip
\item If $ \beta = 1 $, let $ \gammab(a) \sim \Gammab(a) $ i.e. $ \Prob{ \gammab(a) \in dx} = \Unens{x > 0} x^a e^{-x} \frac{dx}{x} $ for $ a > 0 $, and $ B_{\pm 1} \sim \Ber_{\ensemble{\pm 1} }(1/2) $ independent of $ \gammab(a) $; then,
\begin{align*}
\frac{1}{n^{3/4}} \, M_n^{(1)} \cvlaw{n}{+\infty} \Fb_{\!\! 0 \vphantom{a_a} } :\eqlaw B_{ \pm 1} \, \gammab(1/4)^{1/4} \sim e^{- \frac{x^4}{12 }} \frac{dx}{\Ze_0} .
\end{align*}

\medskip
\item If $ \beta = 1 - \frac{\gamma}{\sqrt{n} } $ with $ \gamma \in \Rr $ fixed, let $ \Fb_{\!\! \gamma} \sim  e^{ - \gamma \frac{x^2}{2} - \frac{x^4}{12} } \frac{dx}{\Ze_\gamma } $; then,
\begin{align*}
\frac{1}{n^{3/4}} \, M_n^{(1 - \gamma / \sqrt{n} )} \cvlaw{n}{+\infty} \Fb_{\!\! \gamma}.
\end{align*}

\medskip
\item If $ \beta > 1 $, let $ \Bd \sim \Ber_{\ensemble{\pm 1} }(1/2) $; then,
\begin{align*}
\frac{1}{n} \, M_n^{(\beta )} \cvlaw{n}{+\infty} t_\beta\, \Bd, \qquad \mbox{where } t_\beta = \tanh(\beta t_\beta).
\end{align*}
\end{enumerate}
\end{theorem}
\end{shaded}

We will moreover prove some new theorems, for instance:
\begin{shaded}
\begin{theorem}[Theorem~\ref{Theorem:MagnetisationBeta>1:Functional} in the sequel]\label{Theorem:Intro:MagnetisationBeta>1:Functional}
If $ \beta > 1 $, setting $\Tb_{\! n}^{(\beta)} := 2 \widetilde{V}_{n, \beta} - 1 $ with $ \widetilde{V}_{n, \beta} \sim \widetilde{\nu}_{n, \beta} $, one has
\begin{align}\label{Eq:Intro:Beta>1:Couple}
\prth{ \sqrt{n} \prth{ \frac{ M_\pe{nt}^{(\beta)} }{ n } - \Tb_{\! n}^{(\beta)} }\! , \, \Tb_{\! n}^{(\beta)} }_{\!\! t \geq 0} 
\cvlaw{n}{+\infty} 
\prth{ \frac{1 + \Bd}{2} G_\beta^{(+)}(t) +  \frac{1 - \Bd}{2} G_\beta^{(-)}(t) , t_\beta \, \Bd}_{\!\! t \geq 0}
\end{align}
where $ \prth{ G_\beta^{(+)}(t) , G_\beta^{(-)}(t) }_{t \geq 0} $ is a 2-dimensional centered Gaussian process independent of $ \Bd $ with covariance
\begin{align*}
\Esp{G_\beta^{(+)}(t)  G_\beta^{(-)}(s) } = \prth{ \frac{1 + t_\beta}{2} }^2 (t \wedge s)
\end{align*}

More generally, there exists a Gaussian sheet $ (\widehat{\Zeb}(p, t))_{p \in [0, 1], t \geq 0} $ that governs (almost) all Gaussian random variables in the Curie-Weiss model (Lemma~\ref{Lemma:CvLawFuncZ} in the sequel).
\end{theorem}
\end{shaded}

The analogous result for $ \beta = 1 $ can be deduced from results of Papangelou \cite{Papangelou} and Jeon \cite{JeonCW} and the case $ \beta < 1 $ is given in \eqref{Eq:Beta<1:CvMagFunctional} and reads, with a Brownian Motion $ W $ (which is a particular fiber of the Gaussian sheet $ \widehat{\Zeb} $) independent of $ G_\beta \sim \Ns\prth{ 0, \frac{\beta}{1 - \beta} } $
\begin{align}\label{Eq:Intro:Beta<1:CvMagFunctional}
\prth{ \frac{M_\pe{nt}^{(\beta)}}{\sqrt{n}} }_{\!\! t \geq 0} 
\cvlaw{n}{+\infty } 
\prth{2 W_t + tG_\beta}_{t \geq 0} 
\end{align}

We now describe the main novelty of the paper~: the method of proof.

\medskip 
\subsection{The almost sure Laplace inversion approach}\label{SubSec:Intro:Laplace}

\subsubsection{\textbf{The Laplace inversion formula for $ \Unens{x > 0} $}}\label{SubSubSec:Intro:Laplace:LaplaceInversion}

The Laplace inversion formula for the modified Heaviside function $ x \mapsto \Unens{x > 0} + \frac{1}{2}\Unens{x = 0} $ reads as follow (see e.g. \cite[(4) p. 130]{Tenenbaum})~:
\begin{align}\label{Eq:LaplaceInversion}
\Unens{x > 0} + \frac{1}{2} \Unens{x = 0} = \int_{c + i\Rr} e^{xs} \frac{d^*s}{s}, \qquad c > 0, \qquad d^*s := \frac{ds}{2i\pi}
\end{align}

It is a consequence of the more general asymptotic formula when $ T\to +\infty $ \cite[Lem.~1.1 p. 131]{Tenenbaum}
\begin{align*}
\int_{ c + i [-T, T] } e^{ x s } \frac{d^*s}{s} 
             = \begin{cases}
                   O\prth{T\inv e^{xc}/\ln(x) } & \mbox{if } x < 0 \\
                   \frac{1}{2} + O\prth{T\inv } & \mbox{if } x = 0 \\ 
                   1 + O\prth{ T\inv e^{xc}/\ln(x) } & \mbox{if } x > 0 
             \end{cases}
\end{align*}

This formula is at the origin of the Perron inversion in the Number Theory literature \cite[\S~2.1 p. 131]{Tenenbaum}~; it also takes the form of a Mellin inversion when $ e^x = y $ (this is the form given in \cite[(4) p. 130]{Tenenbaum}).

\medskip
\subsubsection{\textbf{The uniform coupling of Bernoulli random variables}}\label{SubSubSec:Intro:Laplace:BernoulliCoupling}

It is well known that
\begin{align}\label{Eq:BernoulliCoupling}
B_p \eqlaw \Unens{ U < p}, \qquad U \sim \Us([0, 1])
\end{align}
and that this equality in law gives a coupling of the random variables $ (B_p)_{p \in [0, 1]} $. 

\medskip

Set
\begin{align*}
(U_k)_{1 \leq k \leq n} & \sim \mathrm{i.i.d.}\Us([0, 1]) , 
\qquad
\widetilde{V}_{n, \beta}  \sim \widetilde{\nu}_{n, \beta}
\end{align*}
and
\begin{align*}
B_k(p) := \Unens{ U_k < p} \sim \Ber_\ensemble{0, 1}(p), 
\qquad
X_k(p) := 2B_k(p) - 1 \sim \Ber_\ensemble{\pm 1}(p)  
\end{align*}


The De Finetti characterisation of $ (X_k^{(\beta)})_{1 \leq k \leq n} $ is then equivalent to
\begin{align*}
\boxed{\prth{ X_k^{(\beta)} }_{1 \leq k \leq n} \eqlaw \prth{ \Unens{ U_k < \widetilde{V}_{n, \beta} } }_{1 \leq k \leq n}} \, , 
\qquad 
\widetilde{V}_{n, \beta} \ \mbox{ independent of } \ (U_k)_{1 \leq k \leq n}
\end{align*}

\medskip
\subsubsection{\textbf{The approach}}\label{SubSubSec:Intro:Laplace:Approach}

The conjunction of \eqref{Eq:LaplaceInversion} and \eqref{Eq:BernoulliCoupling} gives rise to 
\begin{align}\label{Eq:BernoulliWithIntegral}
B_p \eqlaw \int_{c + i\Rr} e^{s (p - U) } \frac{d^*s}{s}, \qquad c > 0
\end{align}
since $ \Unens{U = p} = 0 $ a.s. 
In particular, a sum of i.i.d. Bernoulli random variables is given by
\begin{align*}
S_n(p) & := \sum_{k = 1}^n X_k(p) \\
               & = 2 \int_{c + i \Rr} e^{ sp } \prth{ \sum_{k = 1}^n e^{-s U_k} }\frac{d^*s}{s} - n \\
               & =: 2 \int_{c + i \Rr} e^{ sp } Z_n(s) \frac{d^*s}{s} - n
\end{align*}
where $ (Z_n(s))_{s \in c + i\Rr} $ defines the \textit{i.i.d. random function}
\begin{align}\label{Def:IIDfield}
Z_n(s) := \sum_{k = 1}^n e^{-s U_k}  = \bracket{e^{-s \, \cdot } , \, \sum_{k = 1}^n \delta_{U_k} }
\end{align}
and $ \sum_{k = 1}^n \delta_{U_k} $ is the \textit{i.i.d. field}, $ \bracket{\cdot, \cdot} $ being the duality bracket between distributions in $ \De' $ and $ \Ce^\infty_0 $ functions \cite[ch. I-2, p. 21]{SchwartzDistributions}. Since a sequence of random variables $ (U_k)_k $ defines a random field (or a point process) $ \sum_k \delta_{U_k} $, we will use a slight abuse of language and consider a sequence to be a field from now on.

\medskip

Applying this philosophy to the unnormalised magnetisation \eqref{Def:UnnormalisedMagnetisation}, we obtain 
\begin{align*}
M_n^{(\beta)}    :=  \sum_{k = 1}^n X_k^{(\beta)}  
                  =  \sum_{k = 1}^n \prth{ 2\Unens{ U_k <  \widetilde{V}_{n, \beta} } - 1}   
               & =:  \, S_n\prth{\widetilde{V}_{n, \beta} }  
\end{align*}
namely
\begin{align}\label{Eq:MagnetisationWithFields}
\boxed{M_n^{(\beta)} =  \, 2 \int_{c + i\Rr} e^{s  \widetilde{V}_{n, \beta} } Z_n(s) \frac{d^*s}{s} - n}
\end{align}

This formula will be at the core of our analysis of $ M_n^{(\beta)} $~:
\begin{itemize}

\item since $ Z_n $ is a sum of i.i.d.'s, it will converge (with a proper renormalisation) to a Gaussian function $ (\Ze(s))_{s \in c + i\Rr} $ and the limit will be expressed by means of its Laplace transform 
\begin{align}\label{Det:Zhat}
\widehat{\Ze}(p) := \int_{c + i\Rr} e^{s  p } \Ze(s) \frac{d^*s}{s}
\end{align}
that will appear to be a Brownian Bridge at the origin of (almost)
all Gaussian random variables in the Curie-Weiss model: $ (\widehat{\Ze}(p))_{p \in [0, 1]} \sim BB $.

\medskip
\item the randomisation $ \widetilde{V}_{n, \beta} $ will itself satisfy a CLT after a proper rescaling,  

\medskip
\item the conjunction of these two CLTs will give the full phase diagram of the Curie-Weiss model and will allow to prove theorem~\ref{Theorem:MainIntro} in full generality.

\end{itemize}

\newpage
\subsection{Literature comparison}\label{SubSec:Intro:Comparison}

\subsubsection{\textbf{With the surrogate by exchangeability}}\label{SubSubSec:Intro:Comparison:SurrogateExch}

In the companion paper \cite{BarhoumiButzekEichelsbacher}, the notion of \textit{surrogate by exchangeability} was defined using the CLT approximation of the Bernoulli sum followed by the De Finetti randomisation~:
\begin{align}\label{Eq:CLTsurrogate}
S_n(p)  \stackrel{\Le}{\approx} \mu_n(p) + \sigma_n(p) G
\qquad\Longrightarrow\qquad
M_n^{(\beta)} & = S_n\prth{\widetilde{V}_{n,\beta}} \stackrel{\Le}{\approx} \mu_n\prth{\widetilde{V}_{n,\beta}} + \sigma_n\prth{\widetilde{V}_{n,\beta}} G
\end{align}

The principle of the proofs in \cite{BarhoumiButzekEichelsbacher} was to give a precise meaning to the second approximation in law and to analyse the behaviour of the resulting simple random variable composed of two independent explicit random variables. 

We did not use in \cite{BarhoumiButzekEichelsbacher} the functional relation \eqref{Eq:BernoulliCoupling} for the Bernoulli random variables $ B_k(p) $, as a result, the approximation in distribution was done at the level of random variables, without any stochastic process. This is similar in nature, the decomposition allowing for a decoupling through the functions $ \mu_n : p \mapsto \mu_n(p) $ and $ \sigma_n : p \mapsto \sigma_n(p) $, but the introduction of a functional relation is deeper, as always in Probability Theory where stochastic processes allow to explain at a functional level some relations between random variables that may look like peculiarities; this is the ``bijective proof'' approach described e.g. in \cite[\S~4.3]{BianePitmanYor}. 

\medskip

Finally, there is one last computational and theoretical satisfaction in using the Laplace inversion \eqref{Eq:MagnetisationWithFields} and the stochastic process $ \Ze_n $ in place of the CLT surrogate \eqref{Eq:CLTsurrogate}: 
\begin{quote}
\textit{The pre-limit and the limit share the same structural form}.
\end{quote}
 The survival of a certain structure after rescaling/the invariance in form of the prelimits and the limit is a key conceptual feature in a theory. It is the sign that the underlying structure is the ``right'' one to consider. Other examples of the phenomenon are provided by the first author in Random (Unitary) Matrix Theory (``restricted reproducing kernel approach'', \cite{BarhoumiCUErevisited}) and in Integrable probability (``max-independence structure'' \cite{BarhoumiMaxIndep}).

\medskip
\subsubsection{\textbf{Energy-Entropy competition}}\label{SubSubSec:Intro:Comparison:EntropyEnergy}

The classical way to pass to the limit in the total magnetisation is understood through the lens of the \textit{Energy-Entropy competition} for the Curie-Weiss model (and more generally for random walks in random environment), see e.g. \cite{EllisBook, EllisNewman, EllisNewmanRosen}. The previous approach furnishes an alternative to this paradigm, or \textit{I.I.D.-Randomisation competition}, the energy being akin to the i.i.d. field  (the ``order'' that will converge in distribution according to the Donsker invariance principle towards a functional of the Brownian motion/bridge) and the entropy to the randomisation (the ``disorder''). As explained in the companion paper \cite[\S~1.2.2]{BarhoumiButzekEichelsbacher}, the survival of a certain part of the disorder at the limit is named with the conventions of Statistical Mechanics a \textit{marginally relevant disordered system}. This is the case of the sub-critical case $ \beta < 1 $ after passing to the (rescaled) limit in \eqref{Eq:MagDecomposition}, when the randomisation converges in law to a Gaussian random variable. Note that the form of \eqref{Eq:MagDecomposition} already has the structure of the limit, a sum of two independent Gaussian random variables, but this independence only appears at the limit when the randomisation is ``washed away'' in the Brownian Bridge.

\newpage
\subsection{Extensions of the approach}\label{SubSec:Intro:ExtensionApproach}

\subsubsection{\textbf{Change of spins distribution in the Curie-Weiss model}}\label{SubSubSec:Intro:Extension:ChangeOfDistribution}

We would like to consider ``spins'' that do not have values in $ \{0, 1 \} $, i.e. general random variables $ X_k $ in \eqref{Def:CurieWeiss}. If $ X $ is a random variable with probability cumulative function $ F_X : x \mapsto \Prob{X \leq x} $ and $ U \sim \Us([0, 1]) $, it is well known that
\begin{align*}
 F_X(X) \eqlaw U
\qquad\Longleftrightarrow\qquad
 X \eqlaw F_X\inv(U)
\end{align*}
where $ F_X\inv $ is the inverse cumulative function, defined to be right continuous with our convention for $ F_X $. 

The previous case of $ B_p \sim \Ber_\ensemble{0, 1}(p) $ was of this type as one has $ \Prob{B_p = 1} = p = 1 - \Prob{B_p = 0} $ and the inversion of $ F_p := F_{B_p} $ is defined to be $ F_p\inv :x \longmapsto \Unens{ x > 1 - p} $ 
%
%
%
%
%
%
which implies that $ B_p \eqlaw \Unens{ U > 1 - p} = \Unens{1 - U < p} \eqlaw \Unens{U < p } $. 

As a result, instead of the Bernoulli coupling \eqref{Eq:BernoulliCoupling}, one can consider an i.i.d. sequence $ (X_k)_k $ related to $ (U_k)_k $ by $ X_k = F_X\inv(U_k) $. Writing $ F_X\inv $ as the Laplace transform of its inverse Laplace transform 
\begin{align*}
F_X\inv(u) = \int_{c + i\Rr} e^{-su} L_X(s) \frac{d^*s}{s}, \qquad L_X(s) := \Le F_X\inv(s) := \int_{[0, 1]} e^{-u p} F_X\inv(p) dp
\end{align*}
one can perform the previous steps and get an equivalent of the representation \eqref{Eq:MagnetisationWithFields} once a De Finetti representation is found~; this requires to define the exponential bias $ X_k[p] $ of $ X_k $, i.e. the random variable with law proportional to $ e^{p X_k} \bullet \Pp_{X_k} $ and to adapt the proof given for instance in \cite[Lem. B.1]{BarhoumiButzekEichelsbacher}. The $ \{0, 1 \} $-Bernoulli case corresponds to $ L_{B_p}(s) = e^{sp} $. Since the formalism is exactly the same as the one given in the Bernoulli case with a slight modification (replacing $ e^{s \widetilde{V}_{n, \beta}} $ with $ L_{X[p]}(s)\vert_{p = \widetilde{V}_{n, \beta}} $), we will only focus on the Bernoulli case in this paper. 


\medskip
\subsubsection{\textbf{Replacing the randomisation by a randomisation field}}\label{SubSubSec:Intro:Extension:ChangeOfRandomisation}

Instead of the Curie-Weiss model, consider a sequence of spins $ (X_\lambda^{(\beta)})_{\lambda \in \Lambda} $ on a lattice $ \Lambda $ (for instance $ \intcrochet{1, N}^d $ for $d \geq 1$) satisfying
\begin{align}\label{Def:SpinsWithRandomisationField}
\boxed{X_\lambda^{(\beta)} \eqlaw 2 \Unens{ U_\lambda < \, V_{\lambda, \beta} } - 1 } 
\end{align}
for an i.i.d. sequence $ (U_\lambda)_{\lambda \in \Lambda} $ (or \textit{i.i.d. field}) and a \textit{randomisation field} $ (V_{\lambda, \beta})_{\lambda \in \Lambda} $. Then, the previous approach leads to 
\begin{align*}
M_n^{(\beta)}  & :=  \sum_{\lambda \in \Lambda}  X_\lambda^{(\beta)} \\
                & =  \, 2 \int_{c + i\Rr} \sum_{\lambda \in \Lambda} e^{s  V_{\lambda, \beta} } e^{-s U_\lambda} \frac{d^*s}{s} - \abs{\Lambda} 
\end{align*}

Since $ \sum_\lambda a_\lambda b_\lambda = \bracket{ \sum_\lambda a_\lambda \varphi_\lambda, \sum_\mu b_\mu \varphi_\mu } $ for any family $ (\varphi_\lambda)_\lambda $ orthonormal for a scalar product $ \bracket{\cdot, \cdot} $, i.e. $ \bracket{\varphi_\lambda, \varphi_\mu} = \Unens{\lambda = \mu} $, one still has a \textit{decoupling of fields}, i.e. a formula that allows to analyse separately the i.i.d. field $ \sum_\lambda \delta_{ U_\lambda } $ and the randomisation field $ \sum_\lambda \delta_{  V_{\lambda, \beta} } $. If for instance one chooses a scalar product of the form $ \bracket{f, g} := \int f(\zb)\overline{g(\zb)} \mub(d\zb) $ for a good measure $ \mub $, one obtains
\begin{align*}
M_n^{(\beta)}   = \, 2 \int_{c + i\Rr} \int \prth{  \sum_{\lambda \in \Lambda} e^{s  V_{\lambda, \beta} }\varphi_\lambda(\zb) }  \prth{ \sum_{\nu \in \Lambda} e^{-s \, U_{\nu}}\overline{ \varphi_\nu(\zb) } } \mub(d\zb) \frac{d^*s}{s} - \abs{\Lambda} 
\end{align*}
where the two random functions of $ (s, \zb) $ in the parentheses are independent.

The sum of i.i.d.s will still lead to a limiting Gaussian function in $ (s, \zb) $  after renormalisation, and the analysis can be focused on the ``randomisation process'' $ (\sum_{\lambda \in \Lambda} e^{s  V_{\lambda, \beta} }\varphi_\lambda(\zb))_{s, \zb} $.

\medskip

The property \eqref{Def:SpinsWithRandomisationField} is not a peculiarity. It is for instance present in the Ising model (see Annex~\ref{Sec:Ising} for a proof and \cite{Ising, McCoyIsing, McCoyWu} for references on this model). This paper will be focused on the sole Curie-Weiss model where the randomisation field is ``global'', but we plan to tackle in a sequel the case of ``local'' randomisation fields.

\medskip
\subsection{Organisation of the paper}\label{SubSec:Intro:Organisation}

The plan of the paper is as follows~:
\begin{itemize}

\item In \S~\ref{Sec:Prerequisites}, we analyse the relevant Gaussian analytic function $ \Ze $ at stake and its relatives, notably $ \widehat{\Ze} $ defined in \eqref{Det:Zhat} that will be proven to be a Brownian Bridge. 

\medskip
\item In \S~\ref{Sec:Main}, we use the newly introduced paradigm to prove theorem~\ref{Theorem:MainIntro} that gives the full phase transition of the Curie-Weiss model .

\medskip
\item In the short \S~\ref{Sec:Process}, we treat the case of the Donsker-type (functional) rescaling of $ M_n^{(\beta)} $, restricted to the super-critical case. 

\medskip
\item In Annex~\ref{Sec:Ising}, we define the Ising spins and show that they satisfy \eqref{Def:SpinsWithRandomisationField}. 

\end{itemize}

\medskip\medskip\medskip
\section{Prerequisites}\label{Sec:Prerequisites}

\subsection{A Gaussian analytic function}\label{SubSec:Prerequisites:GAF}

We define
\begin{align}\label{Def:CenteredZn}
\centered{Z}_n(s) := Z_n(s) - \Esp{Z_n(s)}, \qquad \Esp{Z_n(s)} = n \frac{1 - e^{-s}}{s} =: n M_\Ze(s)
\end{align}


\medskip
\begin{shaded}
\begin{lemma}[Convergence in law of $ \centered{Z}_n/\sqrt{n} \in \Ce $]\label{Lemma:CvLawZ}
We have the following convergence of the full process on $ \Ce(\Rr, \Rr) $~:
\begin{align*}
\frac{\centered{Z}_n }{\sqrt{n}}  \cvlaw{n}{+\infty}   \Ze  \sim \Ns(0, C_\Ze)
\end{align*}
with
\begin{align*}
C_\Ze(s, w) := \Esp{ \Ze(s)\Ze(w) } =  \frac{1 - e^{-(s + w) }}{s + w} - \frac{1 - e^{-s }}{ s} \frac{1 - e^{-w }}{w} =: M_\Ze(s + w) - M_\Ze(s)M_\Ze(w)
\end{align*}

Moreover, the process can be extended to $ s \in \Cc $ due to the bi-holomorphic (bi-entire) nature of $ C_\Ze $.  
\end{lemma}
\end{shaded}

\medskip

\begin{proof}
For all $ s \in \Rr $, the convergence towards a Gaussian random variable is granted with the CLT for sums of i.i.d.s with finite variance \cite{PetrovSums}. The full process is moreover jointly Gaussian using the vector-valued CLT for $ (\centered{Z}_n(s_1), \dots, \centered{Z}_n(s_k) )/\sqrt{n} $. The convergence in law of finite-dimensional marginals implies the convergence of the full process on the spaces $ \Ce $ and $  \De $ \cite[ch. 7 p. 80 \& 13 p. 138]{Billingsley}. 

It remains to compute the covariance, which is given by the limit of the covariances. But since one has a sum of i.i.d.s, the covariance will be constant in $n$, hence, for all $n$
\begin{align*}
C_\Ze(s, w) & = \frac{1}{n} \Esp{  \centered{Z}_n(s) \centered{Z}_n(w) } \\
              & = \frac{1}{n} \sum_{k = 1}^n \Esp{ \prth{ e^{-s U_k} - \Esp{ e^{-s U_k} } }  \prth{ e^{-w U_k} - \Esp{ e^{-w U_k} } } } \\
              & = \Esp{ \prth{ e^{-s U } - \Esp{ e^{-s U } } }  \prth{ e^{-w U } - \Esp{ e^{-w U } } } } \\
              & =  \Esp{ e^{-(s + w) U } } - \Esp{ e^{-s U } } \Esp{ e^{-w U } } \\
              & = \frac{1 - e^{-(s + w) }}{s + w} - \frac{1 - e^{-s }}{ s} \frac{1 - e^{-w }}{w}
\end{align*}
which furnishes the result.

Last, due to the holomorphicity of $ s \mapsto \frac{1 - e^{-s}}{s} = \sum_{k \geq 1} \frac{(-1)^{k - 1}}{k!} s^{k - 1} $ on the whole complex plane, it is clear that $ C_\Ze(s, \cdot) $ and $ C_\Ze(\cdot, w) $ are holomorphic on $ \Cc $. As  a result, the process can be extended into a holomorphic process on $ \Cc $ and in particular on any line $ c + i\Rr $. 
\end{proof}

\medskip

\begin{remark}
A representation of the random holomorphic function $ \Ze $ is given by the series  
\begin{align*}
(\Ze(s))_{s\in \Cc} \eqlaw \prth{ \sum_{k \geq 0} \frac{G_k}{k!} \, (-s)^k }_{\!\!\! s\in \Cc}
\end{align*}
where $ (G_k)_k $ is a sequence of correlated real (centered) Gaussian random variables with covariance given by
\begin{align*}
\Esp{G_k G_\ell} := \frac{1}{k + \ell + 1} - \frac{1}{(k+1)(\ell + 1)} 
\end{align*}

Verification is left as an easy exercise. 
\end{remark}

\medskip
\subsection{The Brownian Bridge}\label{SubSec:Prerequisites:BB}

We now analyse the Gaussian stochastic process
\begin{align*}
\prth{ \widehat{\Ze}(t) }_{t \in [0, 1]} := \prth{\int_{c + i\Rr} e^{st} \Ze(s) \frac{d^*s}{s} }_{\!\! t \in [0, 1] }
\end{align*}

It is clear that $ \widehat{\Ze} $ is a centered Gaussian process since $ f \mapsto \int_{c + i\Rr} e^{s\cdot } f(s) \frac{d^*s}{s} $ is a linear form. Moreover, it is clearly real since $ \Ze $ is analytic on the whole complex plane and real on $ \Rr $ (hence, no singularity can be picked when moving the line of integration). It is well-defined if one shows that $ \Esp{ \widehat{\Ze}(t)^2 } < \infty $. For this, we compute the covariance of $ \widehat{\Ze} $ for all $ s, t \in [0, 1] $~:
\begin{align*}
C_{\widehat{\Ze}}(s, t) & := \Esp{ \widehat{\Ze}(s)\widehat{\Ze}(t) } \\
                 & = \Esp{ \int_{c + i\Rr} e^{zt} \Ze(z) \frac{d^*z}{z} \int_{c + i\Rr} e^{ws} \Ze(w) \frac{d^*w}{w}} \\
                 & = \int_{(c + i\Rr)^2} e^{zt + ws} C_\Ze(z, w) \frac{d^*z}{z} \frac{d^*w}{w} \\
                 & = \int_{(c + i\Rr)^2} e^{zt + ws} M_\Ze(z + w) \frac{d^*z}{z} \frac{d^*w}{w} -   \int_{ c + i\Rr } e^{zt  } M_\Ze(z ) \frac{d^*z}{z} \int_{ c + i\Rr } e^{ws } M_\Ze(w ) \frac{d^*w}{w} \\
                 & = \int_{(c + i\Rr)^2} e^{zt + ws} \Esp{ e^{-(z + w) U} } \frac{d^*z}{z} \frac{d^*w}{w} -  \int_{ c + i\Rr } e^{zt  } \Esp{ e^{-z U} } \frac{d^*z}{z}  
                   \int_{ c + i\Rr } e^{ws  } \Esp{ e^{-w U} } \frac{d^*w}{w}  \\
                 & = \Esp{ \Unens{ U < s} \Unens{ U < t} } - \Esp{ \Unens{U < s}}\Esp{ \Unens{U < t}} \\
                 & = \Prob{U < s\wedge t} - \Prob{U < s} \Prob{U < t} \\
                 & = s\wedge t - st
\end{align*}

This is the well-known covariance of a \textit{Brownian Bridge} (BB) \cite[ch. I, p. 37 \& ex. 3.10 p. 39]{RevuzYor}~: if $ (W_t)_{t \in [0, 1]} $ designates a Brownian motion on the probability space $ (\Omega, \Feb, \Pp) $ endowed with the filtration $ \Feb := (\Fe_t)_t := (\sigma(B_s, 0 \leq s \leq t))_t $, the BB $ X = (X_t)_{t \in [0, 1]} \equiv (X_t^{(0 \to 0)})_{t \in [0, 1]} $ is defined as $ (W \vert W_1 = 0) \eqlaw (W_t - t W_1)_t \eqlaw (t W_{t\inv - 1})_t $ for the filtration $ \hat{\Feb} := \sigma(\Feb \cup \sigma(W_1)) $. 

Such a BB already appeared in the work of Papangelou \cite{Papangelou} for the critical Curie-Weiss model. Coming back into the physical space (after Laplace inversion), it also comes from the Bernoulli coupling \eqref{Eq:BernoulliCoupling} and the Donsker-Skorokhod-Kolmogorov theorem on the functional limit in law of the empirical distribution on the space $ \De $ \cite{Dudley}.
 One goal of this work is to show that this is this very BB that controls the whole model, at criticality and elsewhere, and that the Gaussian random variables that appear are particular values of this process.

\medskip
\subsection{A deformation of the Brownian sheet}\label{SubSec:Prerequisites:BrownianSheet}

We now rescale $ \centered{Z}_n $ in a process manner, i.e. we consider $ \centered{Z}_\pe{nt}(s)/\sqrt{n} $ for $ t \in [0, 1] $. By doing so, we obtain a random field indexed by $ s \in c + i\Rr \subset \Cc $ and $ t $. To do so, one could use the Donsker invariance principle. But for such a simple sum of i.i.d.'s, it is enough to invoke the vector-valued CLT for the finite-dimensional marginals in points $ (t_i, s_i) $, which gives automatically the answer. The covariance of such a limiting process $ (\Zeb_t(z))_{t \in [0, 1], z \in \Cc} $ is given by  
\begin{align*}
C_\Zeb(z, t \vert w, s) & := \Esp{ \Zeb_t(z) \Zeb_s(w) } \\
                 & = \lim_{n \to +\infty} \frac{1}{n} \Esp{ \centered{Z}_\pe{nt}(z) \centered{Z}_\pe{ns}(w) } \\
                 & = \lim_{n \to +\infty} \frac{1}{n} \sum_{k = 1}^{\pe{nt} \wedge \pe{ns}} \Esp{ \prth{ e^{-z U_k} - \Esp{ e^{-z U_k} } }  \prth{ e^{-w U_k} - \Esp{ e^{-w U_k} } } } \\
                 & = \lim_{n \to +\infty} \frac{\pe{nt} \wedge \pe{ns}}{n} 
                          \times \Esp{ \prth{ e^{-z U } - \Esp{ e^{-z U } } }  \prth{ e^{-w U } - \Esp{ e^{-w U } } } } \\
              & =  s\wedge t \times C_\Ze(z, w)  
\end{align*}


We now compute the law of the process in the Fourier-Laplace space which is the ``original'' space in which the Curie-Weiss model is defined, since one has applied the inverse Laplace transform in \eqref{Eq:LaplaceInversion}. Since the obtained multivariate process $ (\widehat{\Zeb}_t(p))_{t, p \in [0, 1]} $ is Gaussian, it is enough to compute its covariance, which yields to~:
\begin{align*}
C_{\widehat{\Zeb}}(p, t \vert q, s) & := \Esp{ \widehat{\Zeb}_t(p) \widehat{\Zeb}_s(q) } \\
                 & =  \int_{(c + i\Rr)^2} e^{zp + wq} C_{\Zeb}(z,t \vert w, q) \frac{d^*z}{z} \frac{d^*w}{w} \\
                 & = s\wedge t \, \prth{ \int_{(c + i\Rr)^2} e^{zp + wq} C_\Ze(z , w) \frac{d^*z}{z} \frac{d^*w}{w}   } \\
                 & = s\wedge t \, \prth{ q\wedge p - qp }
\end{align*}

The process $ (\Wb_{\!\! s, p})_{s, p \in [0, 1]} $ with covariance $ (s, p ; t, q) \mapsto (s\wedge t)(p\wedge q) $ is the celebrated \textit{Brownian Sheet} \cite[ch. 3, ex. 3.11 p. 39]{RevuzYor}. The fact that on the second ``coordinate'' the covariance is the one of the BB gives one of its natural deformations. 

We sumarize these results into the following~:

\medskip
\begin{shaded}
\begin{lemma}[Functional convergence in law of $ \centered{Z}_n/\sqrt{n}$]\label{Lemma:CvLawFuncZ}
We have the following convergence of the full process on $ \Ce([0, 1]\times \Rr, \Rr) $~:
\begin{align*}
\prth{ \frac{\centered{Z}_\pe{nt}(x) }{\sqrt{n}} }_{t \in [0, 1], x\in \Rr}  \cvlaw{n}{+\infty}   \Zeb := \prth{\Zeb_t(x)}_{t, x}  \sim \Ns(0, C_\Zeb)
\end{align*}
with
\begin{align*}
C_\Zeb(t, z \vert s, w) := (t\wedge s) C_\Ze(z, w)
\end{align*}

This last process can be extended to $ x \in \Cc $ due to the bi-holomorphic (bi-entire) nature of $ C_\Ze $. Its inverse Laplace transform in $x$ is the Gaussian sheet $ (\widehat{\Ze}_t(p))_{t, p \in [0, 1]} $ whose covariance is given by
\begin{align*}
C_{\widehat{\Zeb}}(t, p \vert s, q) := (t\wedge s) C_{\widehat{\Ze}}(p, q) = (t\wedge s)(p\wedge q  - pq)
\end{align*}

\end{lemma}
\end{shaded}

\medskip
\subsection{Convergence of integrals}\label{SubSec:Prerequisites:CvIntegrals}

One would like to write
\begin{align*}
\Le\,\mbox{-}\!\!\!\lim_{n \to +\infty} \int_{c + i\Rr} e^{s p} \frac{\centered{Z}_n(s)}{\sqrt{n}} \frac{d^*s}{s} 
               =  \int_{c + i\Rr} e^{s p} \ \ \Le\,\mbox{-}\!\!\!\lim_{n \to +\infty}\! \frac{\centered{Z}_n(s)}{\sqrt{n}} \,\frac{d^*s}{s} 
\end{align*}
but for this to happen, one would require a dominated convergence, or equivalently a coupling of $ \prth{\centered{Z}_n, \Ze} $ on the same probability space. Such a coupling can be obtained with the Komlos-Major-Tusnady theorem \cite{KomlosMajorTusnadyI, KomlosMajorTusnadyII} but a priori only in $ s \in [0, 1] $~; we thus proceed here differently, by coming back to the random variables.

\medskip

\begin{shaded}
\begin{lemma}[Convergence of random integrals]\label{Lemma:ConvergenceIntegrals}
We have, for all $ p \in [0, 1] $
\begin{align*}
\int_{c + i\Rr} e^{s p} \frac{\centered{Z}_n(s)}{\sqrt{n}} \frac{d^*s}{s} 
               \cvlaw{n}{+\infty} \int_{c + i\Rr} e^{s p} \Ze(s) \frac{d^*s}{s} =: \widehat{\Ze}(p)
\end{align*}

Moreover, this last convergence can be extended for the whole process in $p$.
\end{lemma}
\end{shaded}


\medskip
\begin{proof}
Using the Laplace inversion~\eqref{Eq:LaplaceInversion} and the CLT for sums of independent random variables, we easily obtain the convergence for the random variable. The limiting Gaussian random variable is shown to have the desired covariance by a direct computation similar to the one in \S~\ref{SubSec:Prerequisites:BB}.  The convergence for the process can be proven in a similar way by looking at finite-dimensional marginals and using the CLT for vectors. This concludes the proof.
\end{proof}

\medskip\medskip\medskip
\section{The phase transition picture in the Curie-Weiss model}\label{Sec:Main}

\subsection{The case $ \beta < 1 $}\label{SubSec:Main:Beta<1} 

$ $

\begin{shaded}
\begin{theorem}[Fluctuations of the \textit{unnormalised} magnetisation for $ \beta < 1 $]\label{Theorem:MagnetisationBetaSmaller1}
For $ \beta < 1 $, define $  \Zb_{\!\beta} \sim \Ns\prth{ 0, \frac{1}{1 - \beta} } $. Then, one has  
\begin{align}\label{Eq:Beta<1}
\frac{M_n^{(\beta )}}{\sqrt{n}} \cvlaw{n}{+\infty} \Zb_{\!\beta}
\end{align}
\end{theorem}
\end{shaded}


\begin{proof}
Using \eqref{Eq:MagnetisationWithFields}, we have 
\begin{align*}
M_n^{(\beta)}  & = \, 2 \int_{c + i\Rr}  e^{s \widetilde{V}_{n, \beta}} \centered{Z}_n(s) \frac{d^*s}{s} - n + 2 \int_{c + i\Rr}  e^{s \widetilde{V}_{n, \beta}} \Esp{Z_n(s)} \frac{d^*s}{s} \\
                & = \, 2 \int_{c + i\Rr}  e^{s \widetilde{V}_{n, \beta}} \centered{Z}_n(s) \frac{d^*s}{s} + n\prth{    2 \int_{c + i\Rr}  e^{s \widetilde{V}_{n, \beta}} \Esp{e^{-s U }} \frac{d^*s}{s} - 1 }
\end{align*}

We can compute explicitly this last integral with \eqref{Eq:LaplaceInversion}~:
\begin{align*}
\Phi(p) & := 2 \int_{c + i \Rr} e^{ sp} \Esp{ e^{- s U} } \frac{d^*s}{s} - 1 \\
               & \, = \Esp{ 2\Unens{p - U > 0} - 1 } = 2 \Prob{ U < p } - 1 = 2  (p\Unens{0 \leq p \leq 1} + \Unens{p > 1}) - 1 \\
               & \, = (2p - 1) \Unens{0 \leq p \leq 1} + \Unens{p > 1} - \Unens{p < 0} \\
               & \, = 2p - 1
\end{align*}
since $ p \in [0, 1] $.

Define
\begin{align}\label{Def:TnBeta}
\Tb_n^{(\beta)} := 2 \widetilde{V}_{n, \beta} - 1
\end{align}
so that
\begin{align}\label{Eq:MagDecomposition}
\boxed{M_n^{(\beta)} = \, 2 \int_{c + i\Rr}  e^{s \widetilde{V}_{n, \beta}} \centered{Z}_n(s) \frac{d^*s}{s} + n\, \Tb_n^{(\beta)}}
\end{align}

It was proven in \cite[proof thm.~3.1 \& (3.4)]{BarhoumiButzekEichelsbacher} that for $ \beta < 1 $
\begin{align}\label{Eq:beta<1:CvRandomisation}
\sqrt{n}\, \Tb_n^{(\beta)} \cvlaw{n}{+\infty } G_\beta \sim \Ns\prth{0, \frac{\beta}{1 - \beta}}
\end{align}

As a result, 
\begin{align*}
\widetilde{V}_{n, \beta} := \frac{\Tb_n^{(\beta)} + 1}{2} 
                         \stackrel{\Le}{\equivalent{n \to +\infty}} \frac{\frac{G_\beta}{\sqrt{n} } + 1}{2} 
                         \cvlaw{n}{+\infty} \frac{1}{2}
\end{align*}

We now prove that
\begin{align*}
\int_{c + i\Rr}  e^{s \widetilde{V}_{n, \beta}} \frac{\centered{Z}_n(s)}{\sqrt{n}} \frac{d^*s}{s}  
                & \cvlaw{n}{+\infty}  \int_{c + i\Rr}  e^{s / 2} \Ze(s) \frac{d^*s}{s} 
\end{align*}

Using Lemma~\ref{Lemma:CvLawZ} and $ \widetilde{V}_{n, \beta} \to \frac{1}{2} $, one has the convergence in law of the process $ (e^{s \widetilde{V}_{n, \beta}} \frac{\centered{Z}_n(s)}{\sqrt{n}})_s \to (e^{s / 2} \Ze(s))_s $. Moreover, one has
\begin{align*}
e^{-s\widetilde{V}_{n, \beta}} - e^{-s/2} & \ =  -s\prth{\widetilde{V}_{n, \beta} - \frac{1}{2}  } \int_0^1 e^{ -s ( 1/2 + u (\widetilde{V}_{n, \beta} - 1/2) )} du \\
                  & \ = -s\frac{\Tb_n^{(\beta)} }{2} \int_0^1 e^{ -s/2 -us \Tb_n^{(\beta)}/2 } du \\
                  & \stackrel{\Le}{\equivalent{n \to +\infty}}  -s\frac{G_\beta}{2\sqrt{n}} e^{ -s/2  }  
\end{align*}

As a result, one has
\begin{align*}
\int_{c + i\Rr}  e^{s \widetilde{V}_{n, \beta}} \frac{\centered{Z}_n(s)}{\sqrt{n}} \frac{d^*s}{s}
                   & = \int_{c + i\Rr}   \prth{ e^{-s/2} - \frac{G_\beta}{2\sqrt{n}}(1 + o_\Pp(1)) se^{ -s/2  }   } \frac{\centered{Z}_n(s)}{\sqrt{n}} \frac{d^*s}{s} \\
                   & =  \int_{c + i\Rr}  e^{-s/2} \frac{\centered{Z}_n(s)}{\sqrt{n}}  \frac{d^*s}{s}  \prth{ 1 + O_\Pp\prth{ \frac{1}{\sqrt{n} } } } 
\end{align*}
since, differentiating \eqref{Eq:LaplaceInversion} in the sense of distributions, 
\begin{align*}
\int_{c + i\Rr}  e^{-sx} d^*s = -\frac{d}{dx}\int_{c + i\Rr}  e^{-sx} \frac{d^*s}{s} = - \delta_0(dx)
\end{align*}

Using Lemma~\ref{Lemma:ConvergenceIntegrals}, one then obtains
\begin{align*}
\frac{M_n^{(\beta)}}{\sqrt{n}}  & = \, 2 \int_{c + i\Rr}  e^{s \widetilde{V}_{n, \beta}} \frac{\centered{Z}_n(s)}{\sqrt{n}} \frac{d^*s}{s} + \sqrt{n}\,\Tb_n^{(\beta)}  \\
                & \cvlaw{n}{+\infty} \, 2 \int_{c + i\Rr}  e^{s / 2} \Ze(s) \frac{d^*s}{s} + G_\beta \\
                & \qquad\qquad =: 2\widehat{\Ze}(1/2) + G_\beta  
\end{align*}

As $ \widehat{\Ze}(t) \eqlaw t W_{t\inv - 1} $ by construction of the BB \cite[ex. (3.10) p. 39]{RevuzYor}, one has $ 2\widehat{\Ze}(1/2) \eqlaw W_1 \sim \Ns(0, 1) $.

Last, the uniform random variables are independent of the randomisation, hence $ G $  and $ G_\beta $ are independent. This yields
\begin{align*}
G + G_\beta \eqlaw \Zb_\beta \sim \Ns\prth{0, \frac{1}{1 - \beta}}
\end{align*}
which concludes the proof.
\end{proof}

\medskip
\begin{remark}
The decomposition $ \Zb_\beta = G + G_\beta $ shows that part of the randomisation ``field'' survives at the limit and contributes to the final Gaussian random variable. As remarked in \cite[rk.~3.2]{BarhoumiButzekEichelsbacher},
\begin{quote}
In the language of statistical mechanics of phase transitions, when a disorder is present in a statistical system and has a marginal effect, one talks about a \textit{marginally relevant disordered system} (...). [In the context of the KPZ equation or random polymers,] this is analogous to the convergence of the rescaled KPZ equation to the Edwards-Wilkinson universality class which is Gaussian.
\end{quote}

The disorder here is the randomisation $ \widetilde{V}_{n, \beta} $ and the ``order'' is given by the random walk $ Z_n $ or its limiting version $ \Ze $. The interpretation by means of $ Z_n $ allows to fully treat the Curie-Weiss model as a model of random walk in a random environment. It is also reminiscent of the \textit{energy-entropy competition}, the energy being associated to the ``order'' $ Z_n $ and the entropy to the environment $ \widetilde{V}_{n, \beta} $, see \cite{EllisBook, EllisNewman, EllisNewmanRosen}.
\end{remark}

\medskip
\subsection{The case $ \beta_n = 1 + \frac{\gamma}{\sqrt{n}} $, $ \gamma \in \Rr $}\label{SubSec:Main:BetaTrans} 

$ $

\begin{shaded}
\begin{theorem}[Fluctuations of the unnormalised magnetisation for $ \beta_n = 1 - \frac{\gamma}{\sqrt{n}} $, $ \gamma \in \Rr $]\label{Theorem:MagnetisationBeta_n}
Let $ \Fb_{\!\! \gamma} $ be a random variable of law given by
\begin{align*}
\Prob{\Fb_{\!\! \gamma} \in dx } :=  \frac{1}{\Ze_{\Fb_{\! \gamma}}} e^{-\frac{x^4}{12} - \gamma \frac{ x^2}{2}} dx, \qquad  \Ze_{\Fb_{\! \gamma}} := \int_{\Rr}  e^{-\frac{x^4}{12} - \gamma\frac{ x^2}{2}} dx
\end{align*}

Then, 
\begin{align}\label{Eq:Beta_n}
\frac{M_n^{(\beta_n )}}{n^{3/4}} \cvlaw{n}{+\infty} \Fb_{\!\! \gamma}
\end{align}

One has moreover with $ \Tb_n^{(\beta)} $ defined in \eqref{Def:TnBeta}, 
\begin{align}\label{Eq:Beta_n:Couple}
\prth{ n^{1/4}\prth{ \frac{M_n^{(\beta_n )}}{n^{3/4}} - n^{1/4}\Tb_n^{(\beta_n)} }, n^{1/4}\Tb_n^{(\beta_n)}}   \cvlaw{n}{+\infty}  (G, \Fb_{\!\! \gamma}) 
\end{align}
where $ G := 2 \widehat{\Ze}(1/2) \sim \Ns(0, 1) $ is independent of $ \Fb_{\!\! \gamma} $. This last result when $ \gamma = 0 $ is related to \cite[Lem.~3.3 p. 125]{EllisNewman} and is equal to \cite[thm.~1 p. 270, case $ q = 1 $]{Papangelou}.
\end{theorem}  
\end{shaded}

\begin{proof}
It was proven by a direct analysis in \cite[proof Thm.~3.3, (3.10) \& (3.17)]{BarhoumiButzekEichelsbacher} that
\begin{align*}
\Fb_{\!\! n, \gamma}  & :=  n^{1/4} \, \Tb_{\! n}^{(\beta_n)} 
                     \cvlaw{n}{+\infty} \Fb_{\!\! \gamma}
\end{align*}

As a result, using \eqref{Eq:MagDecomposition}
\begin{align*}
\frac{M_n^{(\beta_n)}}{n^{3/4}} & = \, \frac{2}{n^{3/4}} \int_{c + i\Rr}  e^{s \widetilde{V}_{n, \beta_n}} \centered{Z}_n(s) \frac{d^*s}{s} + n^{1/4} \, \Tb_n^{(\beta_n)} \\
                    & = \frac{2}{n^{1/4}} \int_{c + i\Rr}  e^{s \, \prth{ 1 + \Tb_n^{(\beta_n)} } /2 } \, \frac{\centered{Z}_n(s)}{\sqrt{n}} \frac{d^*s}{s} + n^{1/4} \, \Tb_n^{(\beta_n)} \\
                    & \stackrel{\Le}{\equivalent{n \to +\infty}} \, \frac{2}{n^{1/4}} \int_{c + i\Rr}  e^{s /2 } \Ze(s) \frac{d^*s}{s} + \Fb_{\!\! \gamma} \\
                    & \qquad\qquad =: \frac{G}{n^{1/4}} + \Fb_{\!\! \gamma}  
\end{align*}
where $ G := 2\,\int_{c + i\Rr}  e^{s /2 } \Ze(s) \frac{d^*s}{s} = 2\widehat{\Ze}(1/2) $ was shown in the proof of theorem~\ref{Theorem:MagnetisationBetaSmaller1} to be $ \Ns(0, 1) $ and independent of the randomisation $ \widetilde{V}_{n, \beta} $ (hence of its limit $ \Fb_{\!\! \gamma} $).

This concludes the proof.
\end{proof}

\medskip

Due to its importance in the literature, we single out the case $ \gamma = 0 $ which corresponds to $ \beta = 1 $ into the following~:

\medskip

\begin{shaded}
\begin{corollary}[Fluctuations of the unnormalised magnetisation for $ \beta = 1 $]\label{Corollary:MagnetisationBeta=1}
Define $ \Fb := \Fb_{\!\!\gamma}\vert_{\gamma = 0} $, i.e. 
\begin{align*}
\Prob{\Fb \in dx } :=  \frac{1}{\Ze_\Fb} e^{-\frac{x^4}{12}} dx, \qquad
\Ze_\Fb := \int_{\Rr} e^{-\frac{x^4}{12}} dx = 3^{1/4} 2^{-1/2} \Gamma(1/4)
\end{align*}

Equivalently, if $ \gammab(a) \sim \Gammab(a) $ for $ a > 0 $ i.e. $ \Prob{ \gammab(a)\in dx} = \Unens{x > 0} x^{a - 1} e^{-x} dx $, and  if $ B_{\pm 1}(1/2) \sim \Ber_{\ensemble{\pm 1} }(1/2) $ is independent of $ \gammab(a) $, one has
\begin{align*}
\Fb \eqlaw B_{ \pm 1}(1/2) \, \gammab(1/4)^{1/4}  
\end{align*}

Then, 
\begin{align*}
\frac{M_n^{(1)}}{n^{3/4}} & \cvlaw{n}{+\infty} \Fb \\
\prth{ n^{1/4}\prth{ \frac{M_n^{(1)}}{n^{3/4}} - n^{1/4}\Tb_n^{(1)} }, n^{1/4}\Tb_n^{(1)}}   &\cvlaw{n}{+\infty}  (G, \Fb), \qquad G \ \mbox{independent of} \ \Fb
\end{align*}
\end{corollary}
\end{shaded}

This last convergence in distribution corresponds to \cite[thm.~1 p. 270, case $ q = 1 $]{Papangelou} and is related to \cite[Lem.~3.3 p. 125]{EllisNewman}.

\medskip
\subsection{The case $ \beta > 1 $}\label{SubSec:Main:Beta>1}

We recall that the transcendent equation
\begin{align}\label{Def:CriticalPointEquation}
\tanh(x) = \frac{x}{\beta}, \qquad \beta > 1
\end{align}
has two solutions $ \pm x_\beta $ with $ x_\beta > 0 $ (and in fact $ x_\beta \in (1, \beta) $).

We define
\begin{align*}
\Bd \sim \Ber_{\pm 1}\prth{\tfrac{1}{2}}, \qquad 
t_\beta := \frac{x_\beta}{\beta} = \tanh(x_\beta) \in (0, 1] 
\end{align*}

\medskip
\begin{shaded}
\begin{theorem}[Fluctuations of the unnormalised magnetisation for $ \beta > 1 $]\label{Theorem:MagnetisationBeta>1}
If $ \beta > 1 $, one has 
\begin{align}\label{Eq:Beta>1}
\frac{ M_n^{(\beta)} }{ n } \cvlaw{n}{+\infty}  t_\beta\, \Bd
\end{align}

Moreover, 
\begin{align}\label{Eq:Beta>1:Couple}
\prth{ \sqrt{n} \prth{  \frac{ M_n^{(\beta)} }{ n } -  \Tb_{\! n}^{(\beta)} }\! , \, \Tb_{\! n}^{(\beta)} } \cvlaw{n}{+\infty} \prth{ \frac{1 + \Bd}{2} G_\beta^{(+)} +  \frac{1 - \Bd}{2} G_\beta^{(-)}  , t_\beta \, \Bd}
\end{align}
where
\begin{align*}
\prth{ G_\beta^{(+)} , G_\beta^{(-)} } := \prth{ 2\widehat{\Ze}\prth{ \frac{1 + t_\beta}{2}}\! , 2\widehat{\Ze}\prth{ \frac{1 - t_\beta}{2} } }
\end{align*}
is a centered Gaussian vector independent of $ \Bd $ with covariance
\begin{align*}
\Esp{G_\beta^{(+)}  G_\beta^{(-)} } = \prth{ \frac{1 + t_\beta}{2} }^2
\end{align*}
\end{theorem} 
\end{shaded}

\begin{remark}
The convergence of the couple is the equivalent in the case $ \beta > 1 $ of the previous result which is due to Papangelou \cite[thm.~1 p. 270, case $ q = 1 $]{Papangelou}. To the best of the authors' knowledge, it is new.
\end{remark}

\begin{proof}
It was proven by a direct analysis in \cite[proof Thm.~3.8]{BarhoumiButzekEichelsbacher} that
\begin{align*}
\Tb_{\! n}^{(\beta)}   \cvlaw{n}{+\infty}  t_\beta \, \Bd
\end{align*}

As a result, using \eqref{Eq:MagDecomposition}
\begin{align*}
\frac{M_n^{(\beta )}}{n } & = \, \frac{2}{n } \int_{c + i\Rr}  e^{s \widetilde{V}_{n, \beta }} \centered{Z}_n(s) \frac{d^*s}{s} +   \Tb_n^{(\beta_n)} \\
                    & = \frac{2}{ \sqrt{n} } \int_{c + i\Rr}  e^{s \, \prth{ 1 + \Tb_n^{(\beta)} } /2 } \, \frac{\centered{Z}_n(s)}{\sqrt{n}} \frac{d^*s}{s} +  \Tb_n^{(\beta_n)} \\
                    & \hspace{-0.4cm} \stackrel{\Le}{\equivalent{n \to +\infty}} \, \frac{2}{\sqrt{n} } \int_{c + i\Rr}  e^{s (1 + t_\beta \Bd)/2 } \Ze(s) \frac{d^*s}{s} + t_\beta \Bd 
                                 \ \ \ \mbox{ with Lemma~\ref{Lemma:ConvergenceIntegrals}} \\
                    & \ \  =  \frac{2}{\sqrt{n} } \prth{ \Unens{\Bd = 1}  \int_{c + i\Rr}  e^{s (1 + t_\beta )/2 } \Ze(s) \frac{d^*s}{s}    +   \Unens{\Bd = -1}  \int_{c + i\Rr}  e^{s (1 - t_\beta )/2 } \Ze(s) \frac{d^*s}{s}  } + t_\beta \Bd \\
                    & \ \ =: \frac{2}{\sqrt{n} } \prth{ \Unens{\Bd = 1}  \widehat{\Ze}\prth{\tfrac{ 1 + t_\beta}{2} }    +   \Unens{\Bd = -1}  \widehat{\Ze}\prth{\tfrac{ 1 - t_\beta}{2} }  } + t_\beta \Bd
\end{align*}

Remark that $ \Bd = \Unens{\Bd = 1} - \Unens{\Bd = -1} $ and that $ \Unens{\Bd = 1} + \Unens{\Bd = -1} = 1 $, so that $ \Unens{\Bd = \pm 1} = \frac{1}{2}(1 \pm \Bd) $. One concludes the proof with the covariance of the BB.
\end{proof}

\medskip 
\section{Note on the process rescaling}\label{Sec:Process}


In the method of decoupling through (inverse) Laplace transform that we introduce, the functional rescaling of the Curie-Weiss model amounts to the functional rescaling of sums of independent random variables that is thouroughfully studied in the literature \cite{PetrovSums} and which is sumarised in Lemma~\ref{Lemma:CvLawFuncZ}. Typically, the classical Donsker theorem for $ \sum_{k = 1}^\pe{nt} X_k $ applies in a straightforward way. 

As a result, the only relevant information to fully pass to the limit is the coupling in $n$ of the randomisation $ \widetilde{V}_{n, \beta} $. On could try to find a refined coupling structure in the randomisations, allowing to couple the Curie-Weiss model with $n$ spins to the one with $ n + 1 $ spins with an extra randomisation (we were unable to locate in the literature an example of such coupling, but let us mention \cite[lem. B.3, (59)]{BarhoumiButzekEichelsbacher} by the authors where the question of such a coupling is implicitely given through an identity in law involving the total magnetisation with a randomisation using itself the total magnetisation). This goal, although certainly very interesting, goes beyond the scope of this article. We will simply consider here a sequence of independent such random variables, which implies an absence of functional rescaling (or a trivial one, with a random variable multiplied by a deterministic function).

\medskip

We leave to the interested reader the task of writing the functional equivalent of theorem~\ref{Theorem:MainIntro} and give instead the functional version of theorem~\ref{Theorem:MagnetisationBeta>1} which seems to be (arguably) the most interesting one. 

\medskip

\medskip
\begin{shaded}
\begin{theorem}[Functional fluctuations of the unnormalised magnetisation for $ \beta > 1 $]\label{Theorem:MagnetisationBeta>1:Functional}
If $ \beta > 1 $, 
\begin{align}\label{Eq:Beta>1:Couple}
\prth{ \sqrt{n} \prth{  \frac{ M_\pe{nt}^{(\beta)} }{ n } -  \Tb_{\! n}^{(\beta)} }\! , \, \Tb_{\! n}^{(\beta)} }_{\!\! t \geq 0} 
\cvlaw{n}{+\infty} 
\prth{ \frac{1 + \Bd}{2} G_\beta^{(+)}(t) +  \frac{1 - \Bd}{2} G_\beta^{(-)}(t)  , t_\beta \, \Bd}
\end{align}
where
\begin{align*}
\prth{ G_\beta^{(+)}\!(t) , G_\beta^{(-)}\!(t) }_{\! t \geq 0} := \prth{ 2\widehat{\Ze}_t\prth{ \frac{1 + t_\beta}{2}}\! , 2\widehat{\Ze}_t\prth{ \frac{1 - t_\beta}{2} } }_{\! t \geq 0}
\end{align*}
is a vector of Brownian Motions (up to a multiplicative constant) independent of $ \Bd $ with covariance
\begin{align*}
\Esp{G_\beta^{(+)}\!(t)\,  G_\beta^{(-)}\!(s) } = (t \wedge s) \prth{ \frac{1 + t_\beta}{2} }^2
\end{align*}
\end{theorem} 
\end{shaded}

The proof is a direct adaptation of the proof of theorem~\ref{Theorem:MagnetisationBeta>1} using the functional convergence in lemma~\ref{Lemma:CvLawFuncZ}. Such correlated Brownian Motions do not seem to appear in the work of Papangelou \cite{Papangelou} in the critical regime and 
to the best of the authors' knowledge, they constitute a new result. This is not the case of the analogue of theorem~\ref{Theorem:MagnetisationBetaSmaller1} where the functional rescaling of \eqref{Eq:MagDecomposition} gives
\begin{align}\label{Eq:Beta<1:CvMagFunctional}
\prth{ \frac{M_\pe{nt}^{(\beta)}}{\sqrt{n}} }_{\! t \geq 0} 
\cvlaw{n}{+\infty } 
\prth{2 \widehat{\Ze}_t\prth{\tfrac{1}{2}} + t G_\beta}_{t \geq 0}
\end{align}

Here, $ G_\beta $ is related to the Gaussian random variable in \cite[Lem. p. 271]{Papangelou} and \cite[Lem. 3.3]{JeonCW}, and is given by the limit in law of the randomisation $ \Tb_{\! n}^{(\beta)} $ in \eqref{Eq:beta<1:CvRandomisation}. It is independent of the Brownian Motion $ \widehat{\Ze}_{\cdot}\prth{\tfrac{1}{2}} $.


\medskip
\section{Conclusion and perspectives}\label{Sec:Conclusion}

The new paradigm presented in this work in the sole framework of the Curie-Weiss model not only recovers classical results in a unified manner across the subcritical, critical, and supercritical regimes, but also provides new functional limit theorems that might not have all appeared in the literature. The key structural insight is that both the pre-limit and limiting objects share the same common analytic form, pointing to an underlying universality rooted in the Gaussian processes/sheets here defined. The flexibility of this method opens the way for further investigations beyond the non local/global randomisation setting. In particular, the same philosophy applies to the Ising model with a local randomisation field as outlined in the introduction and in the appendix. These developments will be the subject of forthcoming work.


\appendix
\section{The randomisation field in the Ising model}\label{Sec:Ising}

\subsection{Conventions and definitions}\label{SubSec:Ising:Conventions}

\subsubsection{\textbf{General definitions}}\label{SubSubSec:Ising:Conventions:General}

We use the conventions of \cite[Appendix B.]{BarhoumiButzekEichelsbacher}. We define the \textit{logistic function} $ \psi $ and its bijective inverse $ \psi\inv $ by
\begin{align*}
\psi(\alpha) & := \frac{e^\alpha }{e^{\alpha} + e^{-\alpha} } = \frac{1}{1 + e^{-2 \alpha } } = \frac{1}{2}\prth{ 1 + \tanh(\alpha) }, \qquad\alpha \in \Rr \\
\psi\inv(p) & := -\frac{1}{2} \log\prth{ p\inv - 1 } = \Argtanh(2p - 1), \qquad p \in \crochet{0, 1}
\end{align*}

Noting that $ \psi(0) = \frac{1}{2} $, we define
\begin{align}\label{Convention:Bernoulli}
X \equiv X\prth{\tfrac{1}{2}} = X\crochet{0} \sim \Ber_{\pm 1}\prth{\tfrac{1}{2} }, \qquad X(\psi(\alpha)) = X\crochet{\alpha} \sim \Ber_{\pm 1}\prth{ \psi(\alpha) }
\end{align}

It is well known that the exponential bias of a Bernoulli random variable $ X\crochet{\beta} $ is given by \cite{BarhoumiModOmega}
\begin{align}\label{Eq:BernoulliBias}
\frac{e^{\alpha X[\beta] } }{ \Esp{ e^{ \alpha X[\beta] } } } \bullet \Pp_{ X[\beta] } = \Pp_{ X[\alpha + \beta] } 
\end{align}

Define $ S_n\crochet{\alpha} := \sum_{k = 1}^n X_k\crochet{\alpha} $ and $ S_n \equiv S_n\crochet{0} $. The Bernoulli bias \eqref{Eq:BernoulliBias} implies the Binomial bias~:
\begin{align}\label{Eq:BinomialBias}
\frac{e^{\alpha S_n[\beta] } }{ \Esp{ e^{ \alpha S_n[\beta] } } } \bullet \Pp_{ S_n[\beta] }  = \Pp_{ S_n[\alpha + \beta] }  
\end{align}
and one also has the Gaussian bias, with $ G  \sim \Ns(0, 1) $
\begin{align}\label{Eq:GaussianBias}
\frac{e^{\gamma G  } }{ \Esp{ e^{ \gamma G  } } } \bullet \Pp_{ G  } = \Pp_{ G + \gamma } 
\end{align}

\medskip
\subsubsection{\textbf{The Ising model}}\label{SubSubSec:Ising:Decoupling:Ising}

Let $ \Ge := (\Lambda, \Eee) $ be a graph with vertices $ \Lambda \equiv \Lambda_{n, d} $ (for instance $ \intcrochet{1, N}^d $). We define the equivalence relation $ \sim $ by $ i \sim j \quad \Longleftrightarrow\quad {i, j} \in \Eee $. 

We define the adjacency matrix $ C \equiv C_\Ge $ by
\begin{align*}
C[i, j] :=    \Unens{i \sim j}, 
\qquad
C_\beta := \beta C 
\end{align*}

The Ising model was already defined in \eqref{Def:Ising}. It is the sequence of random variables (spins) $ \Bd^{* \beta } := (B^{* \beta }_\lambda)_{\lambda \in \Lambda} $ whose law can equivalently be written as~:
\begin{align}\label{Def:Ising:bis}
\Pp_{ \prth{ \Bd^{* \beta } } } := \frac{e^{ \frac{1}{2} \transp{\Bd} C_\beta \Bd } }{ \Esp{ e^{ \frac{1}{2} \Bd^* C_\beta \Bd  } } } \bullet \Pp_\Bd , \qquad \Bd := (B_\lambda)_{\lambda \in \Lambda} \sim \IID \Ber_{\{\pm 1\}}\prth{\tfrac{1}{2}}
\end{align}
where $ \transp{\Bb} $ is the transpose of the vector $ \Bb $.

\begin{remark}
Since $ \transp{\Bd}\Bd = \sum_{\lambda \in \Lambda} B_\lambda^2 = \abs{\Lambda} $, one can always modify the diagonal of $ C $ by adding $ \gamma I_{\abs{\Lambda}} $ to turn $ C $ into an invertible matrix, the additional constant $ \gamma\times \frac{\beta \abs{\Lambda}}{2} $ obtained in the exponential being absorbed into the renormalisation constant.
\end{remark}

\medskip
\subsection{The randomisation field}\label{SubSec:Ising:Randomisation}

$ $

\begin{shaded}
\begin{lemma}[Randomisation field of the Ising model]\label{Lemma:RandomisationIsing}
Define $ \Vb := (\Vb_{\!\!x})_{x \in \Lambda_{n, d}} $ by
\begin{align}\label{Def:RandomisationFieldsIsing}
\Prob{\Vb \in d\vb} = \mu_n(d\vb) = \frac{1}{\Zg_{n, \beta}} f_{\Zb_\beta}(\vb) \prod_{\lambda\in \Lambda} \cosh(v_\lambda)   dv_\lambda, 
\qquad f_{\Zb_\beta}(\vb) = \frac{e^{- \frac{\beta\inv}{2} \transp{\vb}C\inv \vb} }{\sqrt{\det(2\pi \beta C)}}
\end{align}
where $ \Zg_{n, \beta} := \Esp{ e^{ \frac{1}{2} \Bd^* C_\beta \Bd  } } $. 

Then, one has the randomisation identity~:
\begin{align*}
\prth{ B_\lambda^{* \beta } }_{\lambda \in \Lambda} 
                \eqlaw \Bd [\Vb ] := \prth{ B_{\!\lambda}[\Vb_{\!\!\lambda}] }_{\lambda \in \Lambda } 
\end{align*}

\end{lemma}
\end{shaded}

\medskip

\begin{proof}
We have 
\begin{align*}
\Esp{ g\prth{ \Bd^{(\beta)} } } & := \frac{1}{\Zg_{n, \beta}} \Esp{ e^{ \transp{\Bd} C_\beta \Bd/2 } g\prth{ \Bd } } \\ 
                & = \frac{1}{\Zg_{n, \beta}} \int_{\Rr^\abs{\Lambda}} \Esp{ e^{ \transp{\vb}\Bd  } g\prth{ \Bd  } }  f_{\Zb_\beta}(\vb) d\vb, \qquad \Zb_\beta \sim \Ns(0, C_\beta)  \\
                & = \frac{1}{\Zg_{n, \beta}} \int_{\Rr^\abs{\Lambda}} \Esp{ e^{ \transp{\vb}\Bd  } } \Esp{  g\prth{ B_\lambda \crochet{ v_\lambda }  }_{\lambda \in \Lambda} } f_{\Zb_\beta}(\vb) d\vb \quad \mbox{ with \eqref{Eq:BernoulliBias} } \\ 
                & = \frac{1}{\Zg_{n, \beta}} \int_{\Rr^\abs{\Lambda}} \Esp{  g\prth{ B_\lambda \crochet{ v_\lambda }  }_{\lambda \in \Lambda} } \prod_{\lambda\in \Lambda} \cosh(v_\lambda)  f_{\Zb_\beta}(\vb) d\vb \\
                & = \int_{\Rr^\abs{\Lambda}} \Esp{  g\prth{ \Bd\crochet{\vb} } } \mu_n(d\vb)
\end{align*}
which gives exactly the desired result.
\end{proof}


\begin{remark}
The density \eqref{Def:RandomisationFieldsIsing} is a tensorial change of probability~:
\begin{align*}
\Prob{\Vb \in d\vb}  
                     & = \frac{1}{\Zg_{n, \beta}} f_\Zb(\vb)  \prod_{\lambda\in \Lambda}  d\sinh(v_\lambda) \\
                     & = \frac{1}{\Zg_{n, \beta}} f_\Zb(\Arg\sinh(\xb))  \prod_{\lambda\in \Lambda}   dx_\lambda, \qquad\mbox{where} \quad \Arg\sinh(\xb) := \prth{\Arg\sinh(x_\lambda)}_{\lambda \in \Lambda}.
\end{align*}

In terms of random variables, we have
\begin{align*}
V_\lambda = \sinh(Y_\lambda), \qquad f_\Yb(\xb) = \frac{1}{\Zg_{n, \beta}} f_\Zb(\Arg\sinh(\xb)) 
\end{align*}

The distribution of $ \Yb $ is similar to the Gaussian distribution, but the quadratic form in the exponential in its density uses $ \Arg\sinh(x_\lambda) $. One can nevertheless perform the usual Gaussian computations in these new variables. 
\end{remark}


\begin{remark}
One can also remark that if $ G \sim \Ns(0, 1) $ independent of $ B \sim \Ber_{\pm 1}(1/2) $, the density of $ G + B $ is given by
\begin{align*}
f_{G + B}(x) =   \frac{1}{\sqrt{2\pi e}}  \cosh(x) e^{- \frac{x^2 }{2} } 
\end{align*}

As a result, the density \eqref{Def:RandomisationFieldsIsing} is also the law of the sum of a correlated Gaussian vector and a sum (involving the covariance matrix $ C_\beta $) of (correlated) Bernoulli random variables (namely the Ising spins)~: $ \Vb \eqlaw \Zb_\beta + C_\beta \Bd^{*\beta} $. The proof is an easy computation left to the interested reader.
\end{remark}

Let $ (U_\lambda)_\lambda \sim \IID \Us([0, 1])$, $ (B_\lambda(p))_\lambda = (\Unens{U_\lambda \leq p})_x $ be an i.i.d.~Bernoulli sequence and $ P_\lambda := \psi(V_\lambda) $. 
We have
\begin{align*}
\prth{ B_\lambda^{*\beta} }_{\lambda \in \Lambda} \eqlaw \prth{ B_\lambda(P_{ \lambda}) }_{\lambda \in \Lambda_n}, 
\qquad \Prob{B(p) = 1} = p 
\end{align*}
namely
\begin{shaded}
\begin{align}\label{EqLaw:IsingSpinsAsRandomisedBernoulli}
\prth{ B_\lambda^{*\beta} }_{\lambda \in \Lambda} 
                \eqlaw  \Bd [\Vb ]  
                \eqlaw \prth{ \Unens{U_\lambda \leq \psi(V_{\lambda}) } }_{\lambda \in \Lambda}
                \eqlaw \prth{ \Unens{W_\lambda  \leq V_{\lambda} } }_{\lambda \in \Lambda}
\end{align}
\end{shaded}
$ \hspace{-0.5cm} $ where
\begin{align*}
W := \psi\inv(U) = \Argtanh(2U - 1) 
\end{align*}

\medskip
\begin{remark}
One also has
\begin{align*}
W  \eqlaw \frac{\Gumbel(1) - \Gumbel'(1)}{2}
\end{align*}
with $\Gumbel(1) $ independent of $ \Gumbel'(1) $ are two Gumbel-distributed random variables. We indeed have, using $ \Gumbel(1) \eqlaw \ln(\ee\inv) $ with $ \ee \sim \Exp(1) $~:
\begin{align*}
\Prob{\Gumbel'(1) - \Gumbel(1) \leq x } & = \Prob{ \ln(\ee/\ee') \leq x } = \Prob{ \ee' > e^{-x }\ee } = \Esp{ e^{-e^{-x} \ee} } = \frac{1}{1 + e^{-x}} 
\end{align*}
and
\begin{align*}
\Prob{W \leq x} = \Prob{U \leq \psi(x)} = \psi(x) = \frac{1}{1 + e^{-2x}} = \Prob{\Gumbel(1) - \Gumbel'(1) \leq 2x}
\end{align*}
\end{remark}

\medskip
\section*{Acknowledgements}

The authors thank M. Butzek, O. H\'enard, J. Steif and Y. Velenik for interesting insights and references.

\bibliographystyle{amsplain}

\end{document}